\documentclass[12pt,a4paper,english]{article}
\usepackage{amsthm}
\usepackage{amssymb} %% Para el conjunto vacio.
\usepackage{stmaryrd}
\usepackage{graphicx}
\usepackage{amsmath}
\usepackage{verbatim}
\usepackage{babel}
\usepackage[utf8]{inputenc}
\newlength{\ancho}
\textwidth=\ancho
\renewcommand{\th}{\theta}

\theoremstyle{remark}

\theoremstyle{plain}
\newtheorem{theorem}{Theorem}
\newtheorem{lemma}[theorem]{Lemma}
\newtheorem{claim}[theorem]{Claim}

\theoremstyle{definition}
\newtheorem{definition}[theorem]{Definition}

\usepackage{amsmath}

\newcommand{\func}{\rightarrow}

\newcommand{\iso}{\cong}

\renewcommand{\o}{\vee}

\newcommand{\y}{\wedge}

\renewcommand{\emptyset}{\raisebox{-0.17ex}{\large$\varnothing$}}
\newcommand{\<}{\langle}
\renewcommand{\>}{\rangle}
\newcommand{\dist}{{\bf dist}}  
\newcommand{\join}{{\bf join}}
\newcommand{\ori}{{\bf ori}}
\newcommand{\onto}{{\bf onto}}
\newcommand{\Abs}{{\bf Abs}}
\newcommand{\Modi}{{\bf Mod1}}
\newcommand{\Modii}{{\bf Mod2}}
\newcommand{\exi}{{\bf exi}}  
\newcommand{\puno}{$\mathbf{p_1}$}  
\newcommand{\pdos}{$\mathbf{p_2}$}  
\newcommand{\keywords}[1]{{\renewcommand{\thefootnote}{\relax}\footnotetext{\emph{Keywords:}
    #1}}}
\begin{document}
\title{Factor Congruences in Semilattices}  
\author{Pedro S\'anchez Terraf\thanks{Supported by Conicet.}}
\maketitle  
\begin{abstract}
We characterize factor congruences in semilattices by using 
generalized notions of order ideal and of direct sum of ideals. When
the semilattice has a minimum 
(maximum) element, these generalized ideals turn into ordinary (dual)
ideals.
\bigskip

En este trabajo damos una caracterización de
las congruencias factor en semire\-tí\-culos usando nociones generalizadas
de ideal y suma directa de ideales. Cuando un se\-mi\-retículo tiene
elemento mínimo (máximo), estos ideales generalizados resultan ideales
(duales) ordinarios.
\end{abstract}


\section{Introduction}
Semilattices are ordered structures that admit an algebraic
presentation. Formally, a semilattice may be presented as a partially
ordered set $\<A,\leq\>$ such that for every pair $a,b\in A$ there
exists the supremum or \emph{join} $a\o b$ of the set $\{a,b\}$. Equivalently, $\<A,\o\>$ is
a semilattice if $\o$ is an idempotent commutative semigroup
operation. In this case we call $A$ a \emph{$\o$-semilattice}
(``join-semilattice''). We can also define a partial \emph{meet}
operation $a\y b$ in a $\o$-semilattice that equals the infimum of
$\{a,b\}$ whenever it exists.

The paradigmatic example of a semilattice is given by any family of
sets that is closed under (finite) union. Actually, the free
semilattice on $n$ generators is given by the powerset of
$\{0,\dots,n-1\}$ minus $\emptyset$. A reference for concepts of
general algebra and ideals in semilattices is \cite{Gratzer}.

The aim of this paper is to obtain an inner characterization of direct product
decompositions of semilattices akin to those in classical algebra. To
attain this goal we
represent these decompositions by means of factor congruences. A
\emph{factor congruence} is  the kernel $\{\<a,b\>\in A : \pi(a) = \pi(b)\}$
 of a projection $\pi$ onto a direct factor of $A$. Thus, a direct
 product representation $A\iso A_1 \times A_2$ is determined by the pair of
 \emph{complementary} factor congruences given by the canonical
 projections $\pi_i : A \func A_i$.

We prove that factor congruences and direct representations can be
described by generalized notions of order ideal and of direct sum of
ideals. When the semilattice has a minimum 
(maximum) element, these generalized ideals turn into ordinary (dual)
ideals. 

The main feature of this characterization is that it is completely
analogous to similar definitions in the realm of classical algebra,
as in ring theory. By using the partial $\y$ operation we devise a
definition with the same (quantifier) 
complexity as the reader may encounter in  classical algebra, and we
obtain an equational relation between each element and its direct
summands in a given decomposition. 

By the end of the paper we apply the characterization to the bounded
case and show that  each of the axioms is necessary.
 
\section{Generalized Direct Sums}

\begin{figure}[h]
\begin{center}
\includegraphics[width=25em,keepaspectratio=true]{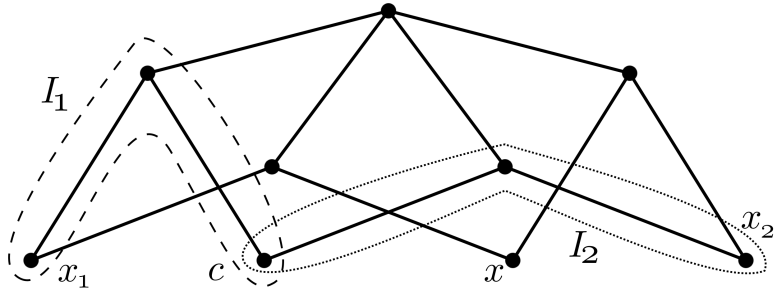}
\end{center}
\caption{A semilattice $A$ isomorphic to $I_1\times I_2$.}\label{fig:paraguas}
\end{figure}

The key idea in our characterization stems from the fact that in a
direct product of join-semilattices there must exist 
``non-trivial'' meets  satisfying certain modularity and absorption laws with respect
to join. And conversely, the existence of meets in a
direct product implies the existence of them in each factor. We state
without proof the latter property: 
\begin{claim} \label{claim: meet_prod} For every $ \o$-semilattices
$C,D$ and elements $c,e \in C$ and $d,f \in D$,     if 
 $ \<c, d\> \y \<e, f\>$ exists in $C \times D$, then $c \y e$ and  $d \y
f$  exist (and conversely).
\end{claim}  
As an a example, a semilattice $A$ is pictured  in
Figure~\ref{fig:paraguas}. $A$ is isomorphic to the direct product of
its subsemilattices $I_1$ and $I_2$. According to this representation,
the pair $\<x_1,x_2\>$ corresponds to $x$. In spite of $x$ not being
expressible in terms of $x_1$, $x_2$ and $\o$, we can recover $x$ as
the infimum  $(x  \o x_1) \y (x \o x_2)$. Other relations between $x$,
$x_1$, $x_2$ and $c$ (the only element in $I_1\cap I_2$) may be found;
next we choose four of them that hold exactly when $x$ corresponds to
$\<x_1,x_2\>$ in a direct decomposition. 

 In the following, we will write formulas in the extended language $
 \{\o, \y \} $. The formula ``$x \y y =  z$'' will be interpreted  as 
``the infimum of $\{x, y\}$ exists and
equals  $z$'':
\[\exists w: \ w\leq x,y \ \ \&\ (\forall u : u\leq x,y\; \rightarrow\; u\leq
w) \ \&\  w = z, \] 
unless explicitly indicated; in general, every equation $t_1 = t_2$
involving $\y$ will be read as ``if either side exists, the other one
also does and they are equal''.  It is easy to see that  the
associative law holds for the partial $\y$ operation in every
$\o$-semilattice, and we will apply it without any further mention.

Heretofore, $A$ will be a $ \o$-semilattice and $c \in A$  arbitrary but
fixed. Let  $\phi_c (x_1, x_2, x) $ be the conjunction of
the following  formulas:
\begin{quote}
\begin{itemize}  
\item [\dist] $  x  = (x  \o x_1) \y (x \o x_2). $
\item [\puno] $  x_1  = (x \o x_1) \y  (c  \o x_1). $ 
\item [\pdos] $  x_2  = (x \o x_2) \y  (c  \o x_2). $ 
\item [\join] $x_1 \o x_2 = x \o c.$  
\end{itemize} 
\end{quote}
\begin{definition}\label{defn:direct_sum}
Assume $I_1, I_2$ are subsemilattices of $A$. We will say that $A$ is
the \emph{$c$-direct sum} of $I_1$ and $I_2$ (notation: $A=I_1 \oplus_c
I_2$) if and only if the following hold: 
\begin{quote}
\begin{itemize}  
\item [\Modi] For all $x,y\in A$, $x_1\in I_1$ and $x_2 \in I_2$,  if
$x\o c \geq x_1\o x_2$ then 
\[ \bigl((x  \o x_1) \y (x \o x_2)\bigr) \o y = (x \o y \o x_1) \y (x
\o y \o x_2).\]
\item [\Modii] For all $x,y\in A$, $x_1\in I_1$ and $x_2 \in I_2$,  if
$x \leq x_1\o x_2$ then 
\[ \bigl((x  \o x_i) \y (c \o x_i)\bigr) \o y = (x \o y \o x_i) \y (c
\o y \o x_i)\]
for $i=1,2$.
\item [\Abs] For all $x_1, y_1 \in I_1$ and $z_2 \in I_2$, we have: $x_1 \y (y_1 \o z_2) = x_1 \y (y_1 \o  c) $ (and
interchanging $I_1$ and $I_2$).
\item [\exi] $ \forall x_1 \in I_1, 
x_2 \in I_2 \; \exists x \in A:  \phi_c (x_1, x_2, x). $ 
\item [\onto]  $ \forall x \in A \ \exists x_1 \in I_1,
x_2 \in I_2 : \; \phi_c (x_1, x_2, x).$  
\end{itemize}
\end{quote}
\end{definition}
In order to make notation    lighter, we will drop the subindex
``$c$''. Let us notice that \exi{} implies: 
\begin{quote}
\begin{itemize}
\item[\ori] $ \forall x_1 \in I_1, x_2 \in I_2 \; (x_1 \o x_2 \geq c)$.
\end{itemize}
\end{quote}
\begin{lemma}\label{l:suma_dir} Assume $A=I_1 \oplus I_2$. Then 
$\phi(x_1,x_2,x)$ defines an isomorphism $\<x_1,x_2\>
\stackrel{\phi}{\mapsto} x$ between $I_1 \times I_2$ and
$A$. 
\end{lemma}
\begin{proof}  
First we'll see that the function $\<x_1, x_2\> \stackrel
 {\phi}{\mapsto}  x $   is well defined.  Let us
 suppose $ \phi (x_1, x_2, x) $ and  $ \phi (x_1, x_2, z) $ (there is  at
 least  one possible image by \exi). We operate as follows:
\begin{align*}  
z &= (z \o x_1) \y    (z \o x_2) && \text{by \dist{} for $z$} \\ 
&= [z \o ((x  \o x_1) \y (c  \o x_1))] \y \\
& \quad \quad \y  [z \o ((x  \o x_2) \y (c  \o x_2))] &&
\text{by \puno{} and  \pdos for $x$.}
\end{align*}
By \join{} for $x$ we may apply \Modii:
\[= [(z \o x \o x_1) \y (\underline{c \o z} \o x_1)] \y 
   [(z \o x \o x_2) \y (\underline{c \o z} \o  x_2)]\]
and by \join{} for $z$,
\[ = [(x \o z \o x_1) \y (x_1 \o x_2)] \y [(x \o z \o x_2) \y (x_1 \o x_2)]  \]
This last term is   symmetric in $x$ and $z$, hence  we
 obtain $x=z$. The function thus  defined by $\phi$ is surjective by \onto; let
 us now prove that it is 1--1. If  $ \phi (x_1, x_2,x )$ and $ \phi (y_1,
 y_2,x) $, we have:  
\begin{align*}  
x_1 &= x_1 \y (x_1 \o x_2) \\ 
&= x_1 \y (x \o c) && \text{by \join} \\ 
&= x_1 \y (y_1 \o y_2) && \text{by \join{} again} \\ 
&= x_1 \y (y_1 \o c) && \text{by \Abs}  
\end{align*}  
hence $x_1 \leq y_1 \o c$. Symmetrically, $y_1 \leq x_1 \o c$ and  in conclusion 
\begin{equation}\label{eq: 1} 
 x_1 \o c = y_1 \o c.
\end{equation}  
On the other hand, 
\begin{align*} 
 x \o y_1 &= ((x  \o x_1) \y (x \o x_2 )) \o y_1 && \text
{by  \dist} \\  
&= (x   \o y_1 \o x_1) \y (x  \o y_1 \o x_2 ) && \text
{by  \Modi} \\  
&= (x \o y_1 \o x_1) \y (x  \o y_1 \o x_2 \o c) && \text{by
\ori} \\ 
 &= (x \o y_1 \o x_1) \y (x \o c) && \text{since $x \o c \geq
x_1,  y_1$} \\ 
 &= x \o y_1 \o x_1 && \text{by the same argument.}
\end{align*}  
Symmetrically,
\begin{equation*}  
x \o x_1= x \o y_1
 \o x_1.
\end{equation*}  
and hence $x \o x_1= x \o y_1$. Collecting this with~(\ref{eq: 1}) and
using \puno{} we have
\[x_1 = (x \o x_1 ) \y (c \o x_1) = (x \o y_1 ) \y (c \o y_1) = y_1.\]
By the same reasoning, $x_2=y_2$.

We will now prove that $\phi$ preserves $\o$. Let us suppose that $ \phi
(x_1, x_2, x)$ and  $ \phi (z_1, z_2, z)$; since each of $I_1, I_2$ is a
subsemilattice,  we know that $x_j
\o z_j \in I_j$ for $j=1,2$. We have to see that $ \phi (x_1 \o z_1, x_2
\o z_2, x \o z)$. The property \join{}
is immediate.
% and then we may apply \Modi{} and \Modii{} for $x\o z$. 
We now prove \dist:
\begin{align}  
x \o z &= \bigl((x \o x_1) \y (x \o x_2)\bigr) \o z  &&  \text
{by \dist{} for $x$} \\ 
&= (x \o x_1 \o z ) \y (x \o x_2 \o z )  &&  \text
{by \Modi} \\ 
&= (x \o z  \o x_1 \o z_1) \y (x \o z  \o x_1 \o z_2) \y \\ 
& \quad \y (x \o z  \o x_2 \o z_1) \y (x \o z  \o x_2 \o
z_2)\label{eq: 5} 
\end{align}
by \dist{} for $z$ and \Modi. Note that 
\begin{align*}  
x \o z  \o x_1 \o z_2 &= x \o c \o z \o c  \o x_1 \o z_2 &&
\text{by \ori, $x_1 \o z_2 = x_1 \o z_2 \o c$} \\ 
&=  {x_1} \o x_2 \o  {z_1} \o z_2 \o  {x
\o z } &&\text{by \join{} two times} \\
 & \geq x \o z  \o x_1 \o
 z_1,  
\end{align*}  
hence we can eliminate two terms  in  equation~(\ref{eq: 5}) and we
obtain \dist{} for $x\o z$ as expected:
\[(x \o z)  = \bigl ((x \o z)  \o (x_1 \o z_1) \bigr) \y \bigl
((x  \o z)  \o (x_2 \o z_2) \bigr). \] 
We may obtain \puno{} and \pdos{}  similarly. We prove \puno:  
\begin{align*}  
x_1 \o z_1  &= ((x   \o x_1) \y (c   \o x_1))\o z_1 &&  \text{by \puno{} for $x$} \\
& = (x \o z_1  \o x_1) \y (c \o z_1  \o x_1) &&  \text{by \Modii} 
\end{align*}
By \puno{} for $z$ followed by \Modii{} in each term of the meet,
\begin{equation}
= ( x \o z \o x_1 \o z_1) \y (x \o c \o x_1  \o z_1) \y  (c \o z 
\o x_1  \o  z_1) \y (c \o x_1 \o z_1 ) \label{eq: 6}
\end{equation}
and this equals
\[ (  ( x \o z) \o (x_1 \o z_1)) \y (c \o (x_1 \o z_1 )).\]
since the last term in~(\ref{eq: 6}) is less than or equal to the  previous ones. 
\end{proof}
Since $\phi(\cdot,\cdot,\cdot)$ defines an isomorphism (relative to
$I_1, I_2$), there are
canonical projections $\pi_j : A \func I_j$ with $j=1,2$ such that
\begin{equation}\label{eq:1}
\forall x\in A, x_1\in I_1, x_2\in I_2 \ : \ \phi(x_1,x_2,x) \iff \pi_1(x) =x_1 \text{ and } \pi_2(x) =
x_2
\end{equation}
holds.
\begin{theorem}\label{th:bijection}
Let $A$ be a $\o$-semilattice and $c\in A$ arbitrary. The mappings 
\[\begin{array}{ccc}
 \<\th, \delta\> & \stackrel{\mathsf{I}}{\longmapsto} & \<I_\th, I_\delta\> \\
 \<\ker \pi_2, \ker \pi_1\> & \stackrel{\mathsf{K}}{\longmapsfrom} & \<I_1, I_2\>
\end{array}\]
where  $I_\th := \{ a\in A : a\, \th\, c\}$, %$I_2 := \{ a\in A : a
					     %\,\th\,c\}$, is a
					     %bijection 
are mutually inverse maps defined
between the set of pairs of complementary factor
congruences of $A$ and the set of pairs of subsemilattices $I_1, I_2$ of $A$ such
that $A=I_1 \oplus_c I_2$.
\end{theorem}

\begin{proof}
We first prove that map $\mathsf{I}$ is well defined, i.e.  $I_\th
\oplus_c I_\delta= A$ for every pair of complementary factor congruences
$\theta, \delta$. It's clear that $I_\th, I_\delta$ are
subsemilattices of $A$, and we know that the mapping $a\mapsto \<a/\th,
a/\delta\>$ is an isomorphism between $A$ and   $ A/\delta \times
 A/\th$. Under this isomorphism,   $I_\th$ corresponds to $\{\<a', c''\>: a'
 \in A/\delta \}$ 
 and $I_\delta$ to $\{\<c', a''\>: a'' \in A/\th \} $, where $c' = c/\delta$ and
 $c'' =c/\th$. Henceforth we will identify $I_\th$ and $I_\delta$ with
 their respective isomorphic images and we will check the  axioms for $I_\th
\oplus_c I_\delta= A$ in $ A/\delta \times  A/\th$.

In order to see
 \Abs, let  $x_1= \<x, c''\>\in I_\th $, $y_1= \<y, c''\>\in I_\th $ and
 $z_2= \<c', z\>\in I_\delta$
 as in the hypothesis and  suppose that   $x_1 \y (y_1 \o  z_2)$ exists. That is to say, 
\[\<x, c''\> \y [\<y, c''\> \o \<c', z\>] = \<x, c''\> \y \<y \o c', c'' \o
 z\>\] 
exists. By Claim~\ref{claim: meet_prod}, we know that $x \y (y \o c') $
 must  exist in $A/\delta$, and in conclusion,  
\begin{multline*}
\<x, c''\> \y [\<y, c''\> \o \<c', z\>] = \<x \y (y \o c'), c'' \y (c''
\o  z)\> = \\ 
= \<x \y (y \o c'), c''\> = \<x, c''\> \y [\<y, c''\> \o \<c', c''\>].    
\end{multline*}
That is to say, the other meet $x_1 \y (y_1 \o c)$ exists  and equals the former. The other
half is analogous.

To check \Modi, we assume $x= \<x', x''\>$, $y= \<y', y''\>$, $x_1=
  \<x'_1, c''\> \in I_\th$ and $x_2 = \<c', x''_2\>\in I_\delta$. Note that  $x\o c \geq x_1\o x_2$ implies $x'\o c' \geq x'\o x_1'$ and $x''\o c''\geq
  x''\o x_2''$, and hence we have $x'\o y'\o c' \geq x'\o y' \o x_1'$ and
  $x''\o y''\o c'' \geq  x''\o y''\o x_2''$. By applying
  Claim~\ref{claim: meet_prod} we obtain:
  \begin{align*}
    \bigl((x  \o x_1) \y (x \o x_2)\bigr) \o y 
    &=  \<((x' \o  x_1')\y(x'\o c'))\o y', ((x''\o c'')\y (x''\o x_2''))\o
    y''\> \\
    &= \<(x' \o x_1')\o y', (x''\o x_2'') \o y''\> \\
    &= \<(x' \o y' \o x_1') \y (x'\o y' \o c'), (x''\o y'' \o x_2'')\>\\
    &= \<(x' \o y' \o x_1') \y (x'\o y'\o c'), (x''\o y'' \o x_2'') \y (x''\o
    y'' \o x_2'')\>\\
    & = (x \o y \o x_1) \y (x \o y \o x_2).
  \end{align*}
The verification of \Modii{} is similar.

Let $x= \<x', x''\>$, $x_1=
 \<x', c''\> $ and $x_2 = \<c', x''\> $. The map $I_\th\times I_\delta
 \func A/\delta \times A/\theta$ given by 
$\<x_1,x_2\> \mapsto x$ is a bijection. Hence, to check \exi{} and \onto{}
 we just have to check that $\phi_c(x_1,x_2,x)$ holds:
%%. This is straightforward.
\begin{quote}
 \begin{itemize}
 \item [\dist]  We have:  
 \begin{align*}  
 x  &= \<x', x''\> \\
 &= \<x' \o x', x''  \o c''\> \y \<x' \o c', x'' \o x''\> \\ 
 &= [\<x', x''\>  \o \<x', c''\>] \y [\<x', x''\>  \o  \<c', x''\>]\\ 
 &= (x  \o x_1) \y (x  \o x_2),  
 \end{align*}  
 which is what we were looking for.
 \item[\puno] We reason as follows
 \begin{align*}
  x_1  &= \<x', c''\> \\
 &= \<x' \o x', x''  \o c''\> \y \<x' \o c', c''  \o c''\> \\ 
   &= [\<x', x''\>  \o \<x', c''\>] \y [\<c', c''\> \o  \<x', c''\>]\\ 
 &= (x  \o x_1) \y (c  \o x_1),  
 \end{align*}  
 \item [\pdos] It's totally analogous.
 \item[\join] Obvious.  
 \end{itemize}  
\end{quote}
Let $\pi_\th$ and $\pi_\delta$ be the projections from
$A/\delta \times A/\th$ onto $I_\th$ and $I_\delta$, respectively, as
defined by~(\ref{eq:1}). Since $\phi_c( \<x', c''\>,\<c', x''\>,\<x',
x''\>)$, we have 
\begin{equation}\label{eq:2}
\pi_\th(\<x', x''\>) =   \<x', c''\> \text{ and }\pi_\delta(\<x', x''\>) = \<c', x''\>
\end{equation}
for all $x'\in A/\delta$ and $x''\in A/\th$. 

We will prove now that each one of the maps $\mathsf{I}$, $\mathsf{K}$
is the inverse of the other.

%%% 
%%% Proof that ker_delta = theta
We first prove $\mathsf{K}\circ\mathsf{I} = \mathit{Id}$. For this we have to show that
$\ker\pi_\delta = \th$ and $\ker\pi_\th = \delta$. Let $a =
\<a',a''\>$ and $b = \<b',b''\>$. By using~(\ref{eq:2}) we have:
\[\ker\pi_\delta =\{\<a,b\> : \pi_\delta (a) =
\pi_\delta(b)\} =\{\<a,b\> : \<c',a''\> = \<c',b''\>\}=\{\<a,b\> : a''= b''\} \]
i.e., $\ker\pi_\delta = \{\<a,b\> : a/\th= b/\th\} = \th$. The proof for
$\ker\pi_\th = \delta$ is entirely analogous.

%%% 
%%% Proof that I_ker pi2 = I_1
Now we show  $\mathsf{I}\circ\mathsf{K} = \mathit{Id}$. We check $I_{\ker\pi_2} = I_1$.
By taking $x=c$ in \onto{} we can deduce $c\in I_1\cap I_2$,
$\phi(c,c,c)$ and hence $\pi_1(c) = \pi_2(c) = c$. By
inspection it is immediate that $\phi(a,c,a)$ for all $a$; in
particular, for $a\in I_1$. Hence we have $\pi_2(a) = c = \pi_2(c)$ for all $a\in
I_1$ and then $I_1\subseteq I_{\ker\pi_2}$. Now take $a\in
I_{\ker\pi_2}$; we have $\pi_2(a) = c$ and then $\phi(a_1,c,a)$ for
some $a_1\in I_1$. By \join{} we have $a_1\o c = a \o c$. Now
considering \puno{}
 and \dist{} we obtain:
\[a_1 = (a \o a_1) \y (c \o a_1) = (a \o a_1) \y (c \o a) = a, \] 
therefore $a\in I_1$, proving  $I_{\ker\pi_2}\subseteq I_1$. We can obtain  
$I_{\ker\pi_1} = I_2$ similarly.
\end{proof}  
These results can be easily extended for the case of semilattices
$\<A,f_1,f_2,\dots\>$ with
\emph{operators}, i.e., operations on $A$ 
preserving $\o$: $f_i(x\o y) = f_ix \o f_iy$. We only have to include in the
definition of $I_1\oplus_c I_2$ some preservation axioms analogous to
\Modi{} and \Modii. In the case of one unary operator, these are:
\begin{enumerate}
\item $I_1, I_2$ are $f$-subalgebras,
\item For all $x\in A$, $x_1\in I_1$ and $x_2 \in I_2$,  if
$x\o c \geq x_1\o x_2$ then 
\[ f \bigl((x  \o x_1) \y (x \o x_2)\bigr) \geq f(x \o x_1) \y f (x \o x_2).\]
\item For all $x\in A$, $x_1\in I_1$ and $x_2 \in I_2$,  if
$x\o c \leq x_1\o x_2$ then 
\[ f \bigl((x  \o x_i) \y (c \o x_i)\bigr) \geq f (x \o x_i) \y f(c \o x_i)\]
for $i=1,2$. 
\end{enumerate}
If $f$ is just a monotone operation, we may also recover the 
Theorem~\ref{th:bijection} by adding:
\begin{quote}
  For all $x\in A$, $x_1\in I_1$ and $x_2 \in I_2$,  if
  $x\o c = x_1\o x_2$ then $fx\o c = fx_1\o fx_2$.
\end{quote}
\subsubsection*{Unique Factorization and Refinement}
It is well-known that semilattices have the property of unique
factorizability. In general, connected posets\footnote{I.e., such that
  the order relation gives rise to a connected graph.}  have the
\emph{refinement property}: every two finite direct decompositions
admit a common refinement. This was proved by Hashimoto \cite{Hashimoto}.

Furthermore, Chen \cite{Chen} showed that 
semilattices have a ``strong form'' of unique factorizability, namely
the \emph{strict refinement property (SRP)} (see, for instance,
\cite{McKenzie,WillardBFC}). This property is equivalent  to the fact that the
factor congruences form a Boolean algebra. By using our representation
of direct decompositions it can be proved that semilattices  have the
SRP. We now translate the lattice operations on congruences to the set of ``direct
summands''  of $A$, $\{I\leq A: \exists J \leq A, I\oplus J =
A\}$. Set intersection of ideals is the meet. To describe the join,
consider  two pairs of ideals $I_1\oplus I_2 = A = J_1\oplus J_2$. Then, there
exists a subsemilattice $I_2\ovee J_2$ of $A$ such that  $(I_1\cap J_1) \oplus
(I_2\ovee J_2) = A$. Namely, if we assume (by using
Theorem~\ref{th:bijection}) that $A=A_1\times A_2$,
$I_1=\{\<x',c''\> : x'\in A_1\}$ and $I_2=\{\<c',x''\> : x''\in A_2\}$, then 
\[I_2\ovee J_2 := \pi_1(J_2)\times A_2 = \pi_1(J_2)\times
\pi_2(I_2).\]
\section{Some Particular Cases}
In the case our semilattice $A$ has a 0 or a 1 ---minimum or maximum
element, resp.--- the characterization of factor congruences by
generalized ideals takes a much simpler
form. By taking $c$ to be the minimum (maximum), it turns out that the
subsemilattices $I_1$ and $I_2$  are order (dual) ideals. 
\begin{theorem}\label{th:caract_0} Let $A$ be a semilattice with 0. Then $A=I_1 \oplus_0 I_2$  if and
only if $I_1,  I_2 \leq A$ satisfy:   
\begin{quote}
\begin{enumerate}  
\item[\Abs] For all $x_1, y_1 \in I_1$ and $z_2 \in I_2$, we have:
  $x_1 \y (y_1 \o z_2) = x_1 \y y_1  $ (and interchanging $I_1$ and $I_2$).  
\item[\onto]  $ I_1 \o I_2 = A.$  
\end{enumerate}  
\end{quote}
Moreover, $I_1$ and $I_2$ are ideals of $A$.
\end{theorem}  
\begin{proof}
%First, we see that \ori{} holds trivially since $c$ is the minimum. 
%
Since $y_1 \o 0 = y_1$, \Abs{} is equivalent to its new
form as stated; and \Modi, \Modii{} are trivially true. The \join{}
part of $\phi_0$ reduces to  ``$x_1 \o x_2 = x$'' and from
here we conclude that \dist, \puno{} and \pdos{}  hold trivially. Hence
\exi{}  simplifies to
\[\forall x_1 \in I_1, x_2 \in I_2 \; \exists x \in A:\ x_1 \o x_2 =
x,\]
which holds trivially, and \onto{} reduces to:
\[\forall x \in A \exists x_1 \in I_1, x_2 \in I_2 \;  :\ x_1 \o x_2 = x.  \]
and this is equivalent  to $ I_1 \o I_2 = A.$

To see that both of $I_1$ and $I_2$ are ideals, we observe that by Theorem~\ref{th:bijection} each
of them is the $\th$-class of $0$ for some congruence $\th$, hence
they are downward closed. 
\end{proof}
\begin{theorem}\label{th:caract_1}  Let $A$ be a semilattice with
1. Then $A=I_1 \oplus_1 I_2$  if and
only if $I_1,  I_2 \leq A$ satisfy:   
\begin{quote}
\begin{enumerate}  
\item[{\bf exi'}] $ I_1 \o I_2 = \{1\} $  
\item[{\bf onto'}] $I_1 \y I_2 = A$. 
\item[{\bf Mod1'}] $ \forall x_1 \in I_1,
x_2 \in I_2 :\  \forall y\ (x_1\y x_2)  \o y = (y \o x_1) \y
(y \o x_2), $ 
\end{enumerate}  
\end{quote}
Moreover, $I_1$ and $I_2$ are dual ideals (viz.\ upward closed sets)  of $A$.
\end{theorem}  
\begin{proof}
We may eliminate \Abs{} and \Modii{} since they are trivial when $c=1$ (this can be
seen by considering \ori{} after
Definition~\ref{defn:direct_sum}).  Formulas \puno{} and \pdos{}   turn into
\[  x_1  = x  \o x_1 \qquad   x_2  = x  \o x_2 \]
which are  equivalent to say $x_1, x_2 \geq x$. Using this, we reduce
\dist\ to:
\[x = x_1 \y  x_2,\]
and then we  may 
eliminate \puno{} and \pdos{} in favor of this last
formula. \join{} is equivalent to  $x_1 \o
x_2 = 1$, and hence \exi{} reduces to the statement ``for all
$x_1\in I_1, x_2\in I_2$,   $x_1 \o x_2 = 1$ and $x_1\y x_2$
exists''. We may state the first part as  $I_1 \o I_2 = \{1\}$, and
may condense the second part with  \onto{}  
writing  $I_1 \y I_2 = A$.

It remains to take care of \Modi. The hypothesis $x\o c \geq x_1\o
x_2$ trivializes when $c=1$. It is obvious that \textbf{Mod1'} implies
\Modi, by taking $y:=x\o y$. And we may obtain the new version under
the hypotheses of the theorem by taking  $x := x_1 \y x_2$ in \Modi. 
\end{proof}
\section{Independence and necessity of the axioms}
We check independence and
  necessity by providing semilattices with a distinguished element $c$
  and a pair of ideals $I_1$ and $I_2$ such that they satisfy every
  axiom except one and it is not the case that $A=I_1 \oplus_c
  I_2$. 

\begin{figure}[h]
\begin{tabular}{ccc}
\includegraphics[height=10.4em,keepaspectratio=true]{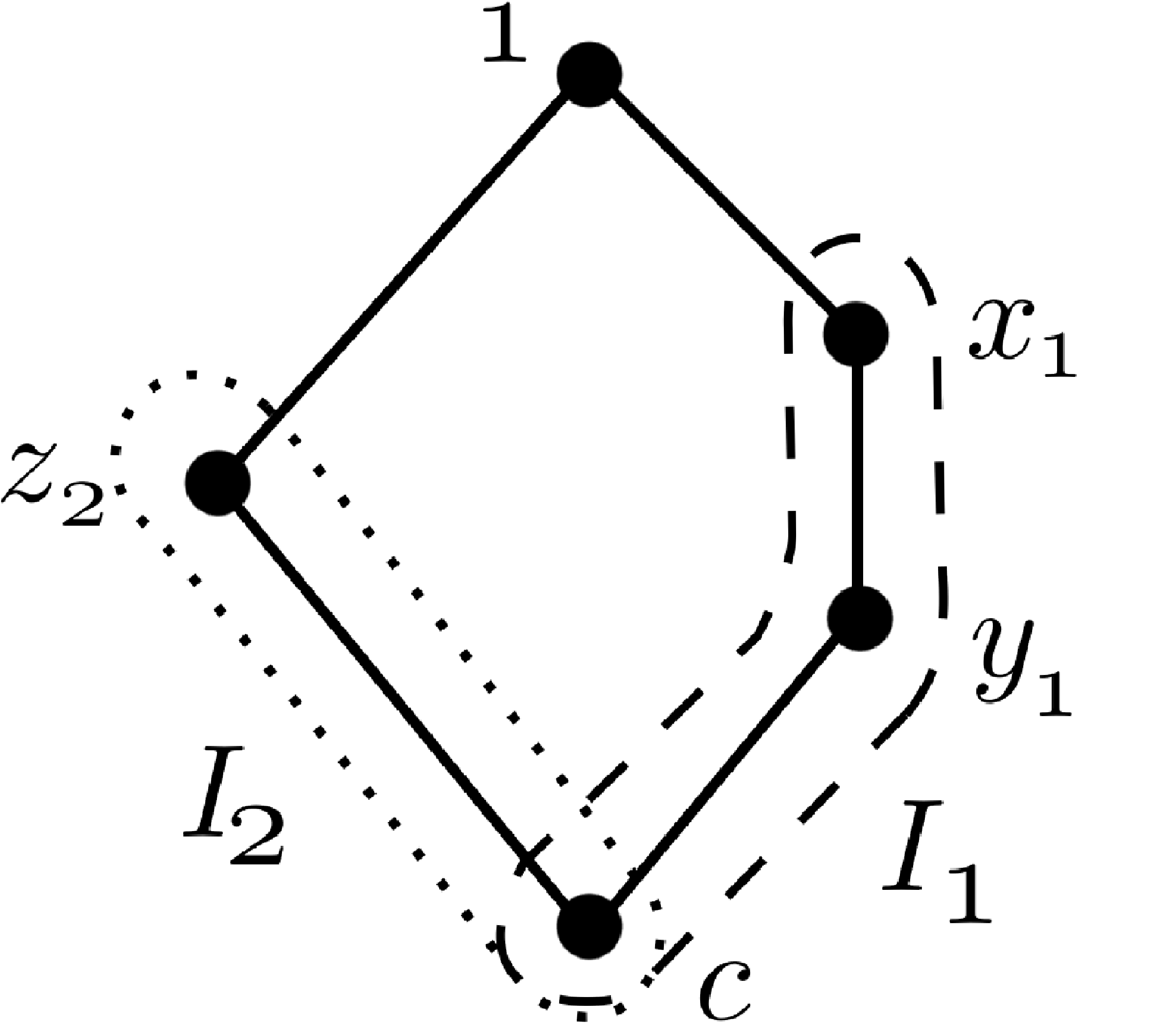} &
\includegraphics[height=10.4em,keepaspectratio=true]{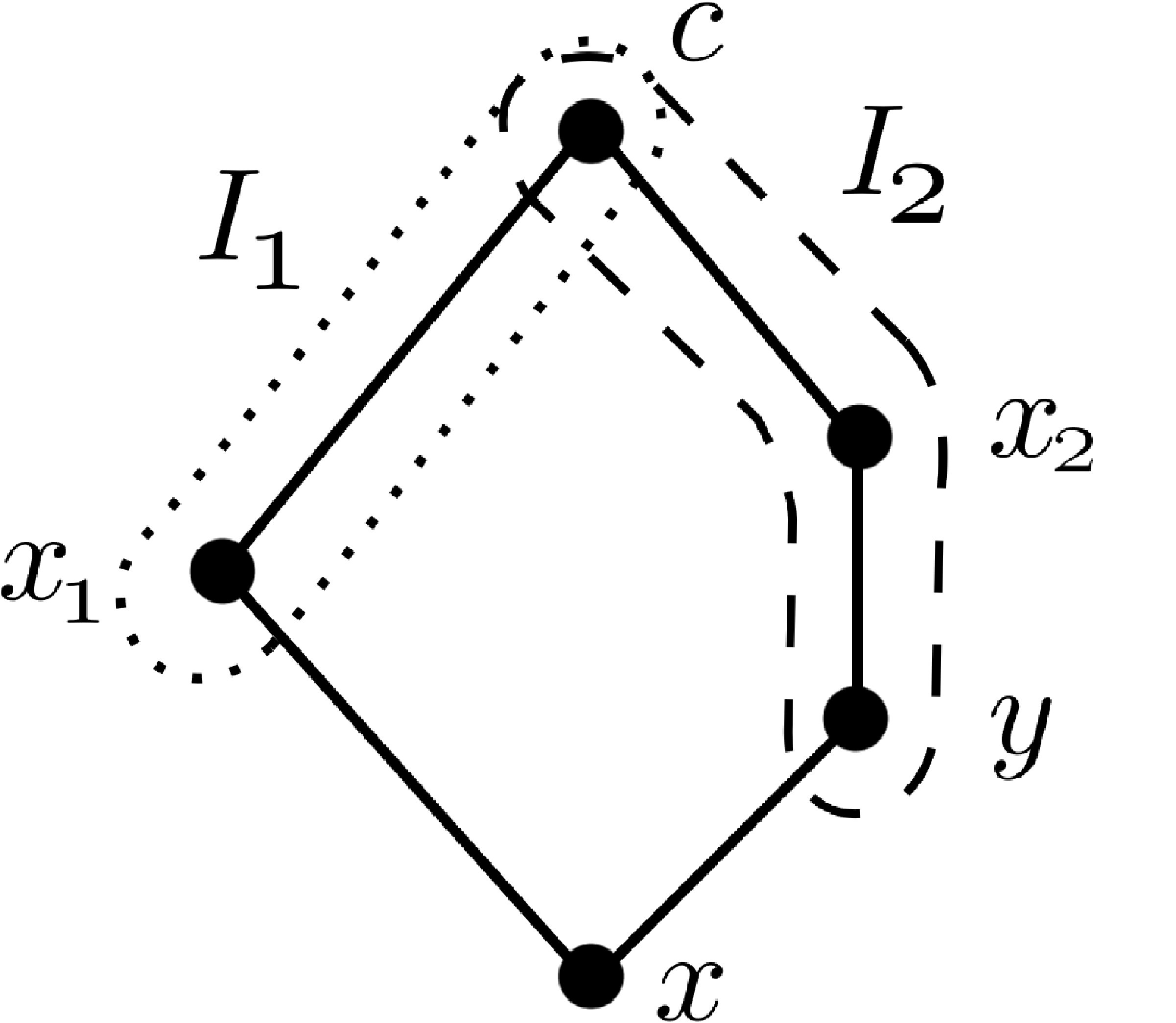} &
\includegraphics[height=10.4em,keepaspectratio=true]{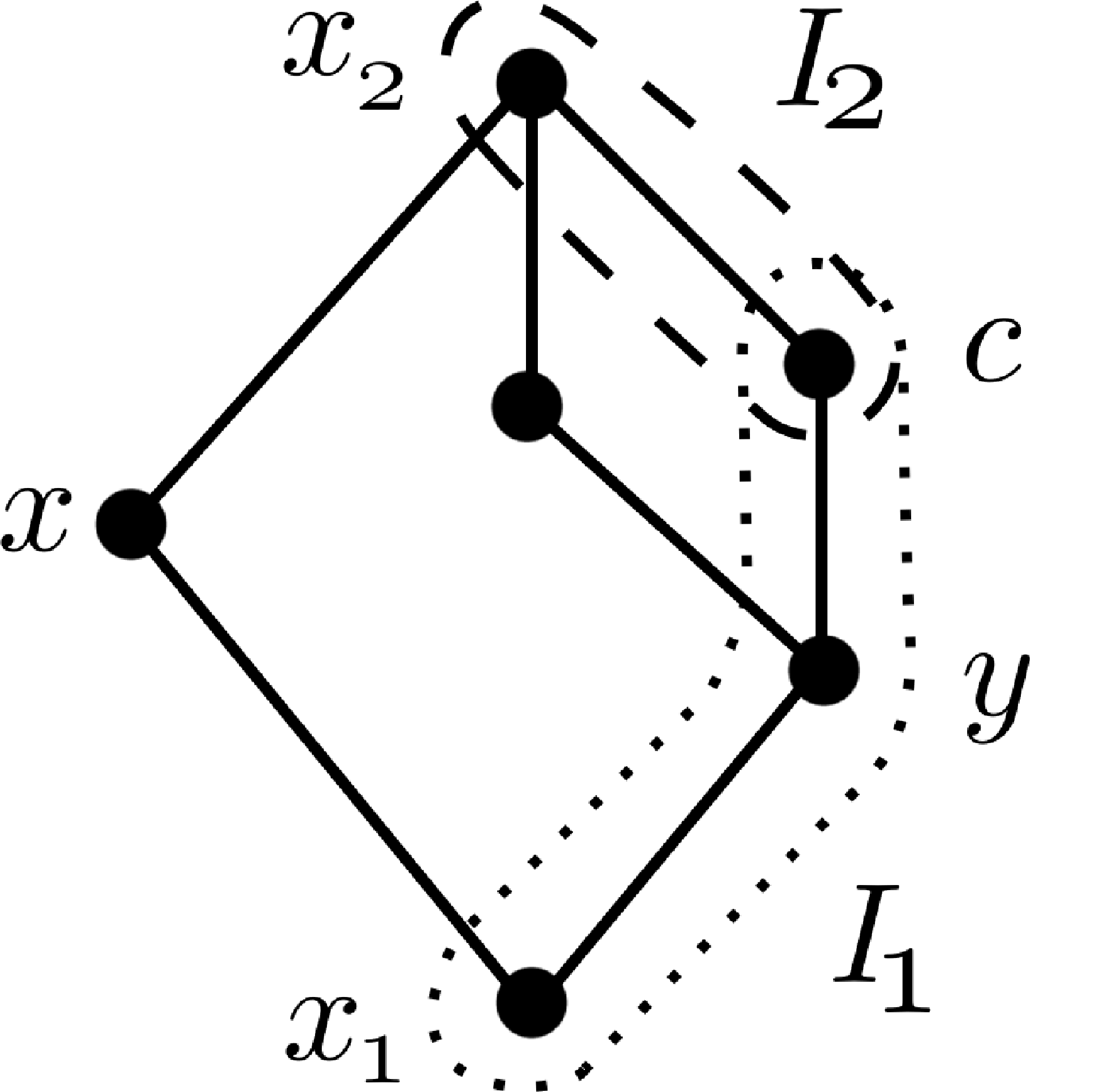}
\\
(a) & (b) & (c)
\end{tabular}
\caption{Counterexamples.}\label{fig:ejemplos}
\end{figure}

It is immediate that the minimal non-modular lattice $N_5$
  provide two such 
counterexamples showing that \Abs{} and \Modi{} are necessary and
independent (see Figure~\ref{fig:ejemplos}~(a) and (b),
respectively). 
In  Figure~\ref{fig:ejemplos}~(a) we have taken $c$ to be the minimum
element, hence Theorem~\ref{th:caract_0} (and its proof) applies and
we only should check \onto{}, which is obvious. Nevertheless, \Abs{}
fails for the labeled elements $x_1$, $y_1$ and $z_2$. Note that
 \Abs{} and \onto{} imply $I_1\cap I_2 = \{c\}$, but this sole assumption does
 not suffices to prove uniqueness of  representations: in
 this example we have $\phi_c(x_1,z_2,1)$ and $\phi_c(y_1,z_2,1)$.

In  Figure~\ref{fig:ejemplos}~(b) $c$ is the maximum of the
semilattice, hence by  Theorem~\ref{th:caract_1} (and its proof) we
only have to check \textbf{exi'} and \textbf{onto'}, again obvious,
and \Modi{} fails for this model.

To see that \Modii{} is necessary and independent is almost as easy as with the other
axioms. Consider the semilattice pictured in Figure~\ref{fig:ejemplos}
(c). First note that every element in $I_2$ is greater than or equal
to every element of $I_1$. From this observation it is immediate to see
that \Modi{} and \Abs{} hold and that $\phi_c(z_1,z_2,z)$ (where
$z_1\in  I_1$ and $z_2\in I_2$) may be simplified as follows:
\begin{align*}
z &= (z \o z_1) \y (z \o z_2) = z \o z_1 & & \text{since } z_1\leq z_2  \\
z_1 &= (z \o z_1) \y (c \o z_1) = z \y c & & \text{by the
    previous line and } c\geq z_1   \\
z_2 &= (z \o z_2) \y (c \o z_2) = (z\o z_2) \y z_2   & & \text{(holds trivially)}  \\
z_1 & \o z_2 = z \o c & &\text{and hence } z_2 = z\o c
\end{align*}
Therefore,  $\phi_c(z_1,z_2,z)$ is equivalent to $(z_1= z\y c) \;\&\; (z_2 = z
\o c)$ and it is now evident that \exi{} and \onto{} also
hold. Finally, consider the labeled elements $x$, $x_1$, $x_2$ and $y$. We
have $x\leq x_1 \o x_2$, and if \Modii{} were true we should have:
\begin{align*}
 \bigl((x  \o x_1) \y (c \o x_1)\bigr) \o y &= (x \o y \o
x_1) \y (c \o y \o x_1) \\
( x \y c ) \o y & = x_2 \y c \\
x_1 \o y &= c\\
y &= c,
\end{align*}
which is an absurdity.

\begin{figure}[h] 
\begin{center}
\begin{tabular}{ccc}
\includegraphics[width=12.5em,keepaspectratio=true]{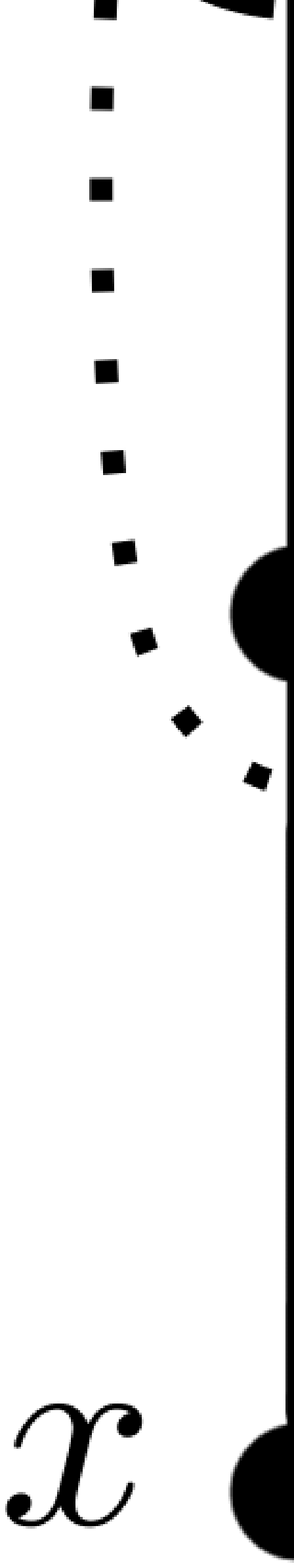} &
\includegraphics[width=12.5em,keepaspectratio=true]{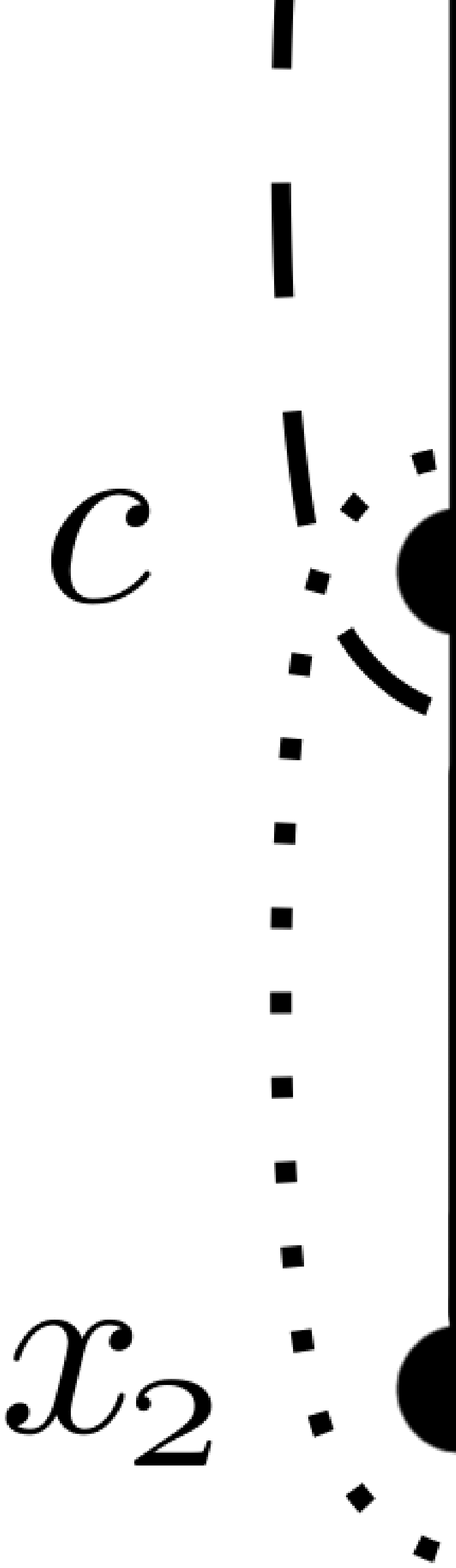} 
\\
(a) & (b) 
\end{tabular}
\end{center}
\caption{Counterexamples.}\label{fig:ejemplos2}
\end{figure}

The last two counterexamples show that
\onto{} and \exi{} are both independent and necessary. In
Figure~\ref{fig:ejemplos2}~(a) there exist no pair $\<x_1,x_2\> \in I_1
\times I_2$ such
that $\phi_c(x_1,x_2,x)$, and in Figure~\ref{fig:ejemplos2}~(b) there
exists no $x$ such that $\phi_c(x_1,x_2,x)$.

\providecommand{\noopsort}[1]{}
\begin{small}\end{small}
\bigskip

\begin{quote}
CIEM --- Facultad de Matem\'atica, Astronom\'{\i}a y F\'{\i}sica 
(Fa.M.A.F.) 

Universidad Nacional de C\'ordoba - Ciudad Universitaria

C\'ordoba 5000. Argentina.

\texttt{sterraf@famaf.unc.edu.ar}
\end{quote}
\end{document}